# Differential calculus on a novel cross-product quantum algebra

**Deepak Parashar**

Department of Mathematics, University of Wales Swansea, Singleton Park, Swansea SA2 8PP, United Kingdom
E-mail: D.Parashar@swansea.ac.uk
URL: http://www-maths.swan.ac.uk/staff/dp/

**Abstract.** We investigate the algebro-geometric structure of a novel two-parameter quantum deformation which exhibits the nature of a semidirect or cross-product algebra built upon $GL(2) \otimes GL(1)$, and is related to several other known examples of quantum groups. Following the R-matrix framework, we construct the $L^{\pm}$ functionals and address the problem of duality for this quantum group. This naturally leads to the construction of a bicovariant differential calculus that depends only on one deformation parameter, respects the cross-product structure and has interesting applications. The corresponding Jordanian and hybrid deformation is also explored.

## 1. The quantum algebra $\mathcal{A}_{r,s}$

The biparametric $q$-deformation $\mathcal{A}_{r,s}$ is defined [1] to be the semidirect or cross-product $GL_r(2) \rtimes_s \mathbb{C}[f, f^{-1}]$ built on the vector space $GL_r(2) \otimes \mathbb{C}[f, f^{-1}]$ where $GL_r(2) = \mathbb{C}[a, b, c, d]$ modulo the relations

$$\begin{array}{ll} ab = r^{-1}ba, & bd = r^{-1}db \\ ac = r^{-1}ca, & cd = r^{-1}dc \\ bc = cb, & [a,d] = (r^{-1} - r)bc \end{array} \quad (1)$$

and $\mathbb{C}[f, f^{-1}]$ has the cross relations

$$\begin{array}{ll} af = fa, & cf = sfc \\ bf = s^{-1}fb, & df = fd \end{array} \quad (2)$$

$\mathcal{A}_{r,s}$ can also be interpreted as a skew Laurent polynomial ring $GL_r[f, f^{-1}; \sigma]$ where $\sigma$ is the automorphism given by the action of element $f$ on $GL_r(2)$. If we let $A = GL_r(2)$ and $H = \mathbb{C}[f, f^{-1}]$, then $A$ is a left $H$-module algebra and the action of $f$ is given by

$$f \triangleright a = a, \qquad f \triangleright b = sb, \qquad f \triangleright c = s^{-1}c, \qquad f \triangleright d = d \quad (3)$$

## 2. The dual algebra $\mathcal{U}_{r,s}$

Knowing properties of cross-product algebras [2, 3], we already know that the algebra dual to $\mathcal{A}_{r,s}$ would be the cross-coproduct coalgebra $\mathcal{U}_{r,s} = U_r(gl(2)) \rtimes_s \mathbb{C}[[\phi]]$ with $\phi$ as an element dual to $f$. As a vector space, the dual is $\mathcal{U}_{r,s} = U_r(gl(2)) \otimes U(u(1))$. Now, the duality relation between $\langle GL_r(2), U_r(gl(2)) \rangle$ is already well-known [4], while that between $\langle \mathbb{C}[f, f^{-1}], U(u(1)) \rangle$ is given by $\langle f, \phi \rangle = 1$, i.e., $U(u(1)) = \mathbb{C}[[\phi]]$. More precisely, we work algebraically with $\mathbb{C}[s^\phi, s^{-\phi}]$ where $\langle f, s^\phi \rangle = s$. This induces duality on the vector space tensor products, the left action dualises to the left coaction, and this results in the dual algebra being a cross-coproduct $\mathcal{U}_{r,s} = U_r(gl(2)) \underset{s}{\ltimes} \mathbb{C}[[\phi]]$. Let us recall that $U_r(gl(2))$, the algebra dual to $GL_r(2)$, is isomorphic to the tensor product $U_r(sl(2)) \otimes \tilde{U}(u(1))$ where $U_r(sl(2))$ has the usual generators $\{H, X_\pm\}$ and $\tilde{U}(u(1)) = \mathbb{C}[[\xi]] = \mathbb{C}[r^\xi, r^{-\xi}]$ with $\xi$ central. Therefore, $\mathcal{U}_{r,s}$ is nothing but $U_r(sl(2))$ and two central generators $\xi$ and $\phi$, where $\xi$ is the generating element of $\tilde{U}(u(1))$ and $\phi$ is the generating element of $U(u(1))$. Also note that $s^\phi$ is dually paired with the element $f$ of $\mathcal{A}_{r,s}$.

## 3. $R$-matrix relations

In the quantum group language, $\mathcal{A}_{r,s}$ is understood as a novel Hopf algebra [5, 6] generated by $\{a, b, c, d, f\}$ arranged in the matrix form

$$T = \begin{pmatrix} f & 0 & 0 \\ 0 & a & b \\ 0 & c & d \end{pmatrix} \quad (4)$$

with the labelling $0, 1, 2$, and $\{r, s\}$ are the two deformation parameters. The $R$-matrix

$$R = \begin{pmatrix} r & 0 & 0 & 0 \\ 0 & \mathsf{S}^{-1} & 0 & 0 \\ 0 & \Lambda & \mathsf{S} & 0 \\ 0 & 0 & 0 & R_r \end{pmatrix} \quad (5)$$

is in block form, i.e., in the order $(00), (01), (02), (10), (20), (11), (12), (21), (22)$ (which is chosen in conjunction with the block form of the $T$-matrix) where

$$R_r = \begin{pmatrix} r & 0 & 0 & 0 \\ 0 & 1 & 0 & 0 \\ 0 & \lambda & 1 & 0 \\ 0 & 0 & 0 & r \end{pmatrix}; \qquad \mathsf{S} = \begin{pmatrix} s & 0 \\ 0 & 1 \end{pmatrix}; \qquad \Lambda = \begin{pmatrix} \lambda & 0 \\ 0 & \lambda \end{pmatrix}; \qquad \lambda = r - r^{-1}$$

The $RTT$ relations $RT_1T_2 = T_2T_1R$ (where $T_1 = T \otimes \mathbf{1}$ and $T_2 = \mathbf{1} \otimes T$) then yield the commutation relations (1) and (2) between the generators. The Hopf algebra structure underlying $\mathcal{A}_{r,s}$ is $\Delta(T) = T \dot{\otimes} T$, $\varepsilon(T) = \mathbf{1}$. The Casimir operator $\delta = ad - r^{-1}bc$ is invertible and the antipode is

$$S(f) = f^{-1}, \quad S(a) = \delta^{-1}d, \quad S(b) = -\delta^{-1}rb, \quad S(c) = -\delta^{-1}r^{-1}c, \quad S(d) = \delta^{-1}a \quad (6)$$

The quantum determinant $\mathcal{D} = \delta f$ is group-like but not central. There are several interesting features of this deformation [5, 6, 1] which cannot be mentioned here. In the $R$-matrix formulation of matrix quantum groups, a basic step is to construct functionals (matrices) $L^+$ and $L^-$ which are dual to the matrix of generators (4) in the fundamental representation. These functions are defined by their value on the matrix of generators $T$

$$\langle (L^\pm)^a_b, T^c_d \rangle = (R^\pm)^{ac}_{bd} \quad (7)$$

where

$$(R^+)^{ac}_{bd} = c^+(R)^{ca}_{db} \quad (8)$$
$$(R^-)^{ac}_{bd} = c^-(R^{-1})^{ac}_{bd} \quad (9)$$

and $c^+$, $c^-$ are free parameters. For $\mathcal{A}_{r,s}$ we make the following ansatz for the $L^{\pm}$ matrices:

$$L^+ = \begin{pmatrix} J & 0 & 0 \\ 0 & M & P \\ 0 & 0 & N \end{pmatrix} \quad \text{and} \quad L^- = \begin{pmatrix} J' & 0 & 0 \\ 0 & M' & 0 \\ 0 & Q & N' \end{pmatrix} \qquad (10)$$

where

$$\begin{aligned}
J &= s^{-\frac{1}{2}(F-A+D-1)} r^{\frac{1}{2}(F-A-D+1)}, & J' &= s^{-\frac{1}{2}(F-A+D-1)} r^{-\frac{1}{2}(F-A-D+1)} \\
M &= s^{-\frac{1}{2}(F-A-D+1)} r^{\frac{1}{2}(-F+A-D+1)}, & M' &= s^{-\frac{1}{2}(F-A-D+1)} r^{-\frac{1}{2}(-F+A-D+1)} \\
N &= s^{-\frac{1}{2}(F+A+D-1)} r^{\frac{1}{2}(-F-A+D+1)}, & N' &= s^{-\frac{1}{2}(F+A+D-1)} r^{-\frac{1}{2}(-F-A+D+1)} \\
P &= \lambda C, & Q &= -\lambda B
\end{aligned}$$

and $\{A, B, C, D, F\}$ is the set of generating elements of the dual algebra $\mathcal{U}_{r,s}$. This is consistent with the action on the generators of $\mathcal{A}_{r,s}$ and gives the correct duality pairings. The commutation algebra is given by the $RLL$ relations $R_{12} L_2^{\pm} L_1^{\pm} = L_1^{\pm} L_2^{\pm} R_{12}$, $R_{12} L_2^+ L_1^- = L_1^- L_2^+ R_{12}$, where $L_1^{\pm} = L^{\pm} \otimes \mathbf{1}$, $L_2^{\pm} = \mathbf{1} \otimes L^{\pm}$, and $R_{12}$ is the same as (5). Finally, we obtain a single-parameter deformation of $U(gl(2)) \otimes U(u(1))$ as an algebra. Including the coproduct, we again obtain $\mathcal{U}_{r,s}$ as a semidirect product $U_r(gl(2)) \underset{s}{\rtimes} U(u(1))$.

## 4. Differential calculus on $\mathcal{A}_{r,s}$

The $R$-matrix procedure [7] is known to provide a natural framework to construct differential calculus on matrix quantum groups. We note here that the $\mathcal{A}_{r,s}$ deformation is not a full matrix quantum group, but an appropriate quotient of one (of multiparameter $q$-deformed $GL(3)$, to be precise). Nevertheless, it turns out that the constructive differential calculus methods [8] work equally well for such quotients. The bimodule $\Gamma$ (space of quantum one-forms $\omega$) is characterised by the commutation relations between $\omega$ and $a \in \mathcal{A}\ (\equiv \mathcal{A}_{r,s})$

$$\omega a = (\mathbf{1} \otimes g) \Delta(a) \omega \qquad (11)$$

and the linear functional $g \in \mathcal{A}'(= \mathrm{Hom}(\mathcal{A}, \mathbb{C}))$ is defined in terms of the $L^{\pm}$ matrices

$$g = S(L^+) L^- \qquad (12)$$

Thus, in terms of components we have

$$\omega_{ij} a = [(\mathbf{1} \otimes S(l_{ki}^+) l_{jl}^-) \Delta(a)] \omega_{kl} \qquad (13)$$

using $L^{\pm} = l_{ij}^{\pm}$ and $\omega = \omega_{ij}$ where $i, j = 1..3$. From these relations, one can obtain the commutation relations of all the left-invariant one-forms with the generating elements of $\mathcal{A}$. The left-invariant vector fields $\chi_{ij}$ on $\mathcal{A}$ are given by the expression

$$\chi_{ij} = S(l_{ik}^+) l_{kj}^- - \delta_{ij} \varepsilon \qquad (14)$$

The vector fields act on the generating elements as

$$\chi_{ij} a = (S(l_{ik}^+) l_{kj}^- - \delta_{ij} \varepsilon) a \qquad (15)$$

Furthermore, using the formula $\mathbf{d}a = \sum_i (\chi_i * a)\omega^i$, we obtain the action of the exterior derivatives ($\mathbf{d}: \mathcal{A} \longrightarrow \Gamma$):

$$\mathbf{d}a = (r^{-2} - 1)a\omega^1 - \lambda b\omega^+ \tag{16}$$
$$\mathbf{d}b = \lambda^2 b\omega^1 - \lambda a\omega^- + (r^{-2} - 1)b\omega^2 \tag{17}$$
$$\mathbf{d}c = (r^{-2} - 1)c\omega^1 - \lambda d\omega^+ \tag{18}$$
$$\mathbf{d}d = \lambda^2 d\omega^1 - \lambda c\omega^- + (r^{-2} - 1)d\omega^2 \tag{19}$$
$$\mathbf{d}f = (r^{-2} - 1)f\omega^0 \tag{20}$$

where $\omega^0 = \omega_{11}, \omega^1 = \omega_{22}, \omega^+ = \omega_{23}, \omega^- = \omega_{32}, \omega^2 = \omega_{33}$. $\mathbf{d}\mathcal{A}$ generates $\Gamma$ as a left $\mathcal{A}$-module, and this defines a first-order differential calclulus $(\Gamma, \mathbf{d})$ on $\mathcal{A}_{r,s}$. The calculus is bicovariant due to the coexistence of the left ($\Delta_L : \Gamma \longrightarrow \mathcal{A} \otimes \Gamma$) and the right ($\Delta_R : \Gamma \longrightarrow \Gamma \otimes \mathcal{A}$) actions. Curiosly, using the Leibniz rule it can be checked that

$$\mathbf{d}(af - fa) = 0, \quad \mathbf{d}(cf - sfc) = 0, \quad \mathbf{d}(bf - s^{-1}fb) = 0, \quad \mathbf{d}(df - fd) = 0, \tag{21}$$

which is consistent with cross relations (2), and so the differential calculus also respects the cross-product structure of $\mathcal{A}_{rs}$.

## 5. The Jordanian deformation $\mathcal{A}_{m,k}$

$\mathcal{A}_{r,s}$ can be contracted [6] (by means of singular limit of similarity transformations) to obtain a nonstandard or Jordanian analogue, say $\mathcal{A}_{m,k}$, with deformation parameters $\{m, k\}$ and the associated $R$-matrix is triangular. In analogy with $\mathcal{A}_{r,s}$, $\mathcal{A}_{m,k}$ can also be considered as the semidirect or cross-product $GL_m(2) \rtimes_k \mathbb{C}[f, f^{-1}]$ built upon the vector space $GL_m(2) \otimes \mathbb{C}[f, f^{-1}]$, where $GL_m(2)$ is itself a Jordanian deformation of $GL(2)$. Thus, $\mathcal{A}_{m,k}$ can also be interpreted as a skew Laurent polynomial ring $GL_m[f, f^{-1}; \sigma]$ where $\sigma$ is the automorphism given by the action of element $f$ on $GL_m(2)$.

## 6. Conclusions

The $\mathcal{A}_{r,s}$ and $\mathcal{A}_{m,k}$ deformations provide interesting new examples of cross-product quantum algebras, both of which have $GL(2) \otimes GL(1)$ as their classical limits. The differential calculus on $\mathcal{A}_{r,s}$ also has an inherent cross-product structure, embeds the calculus on $GL_q(2)$ and is also related to the calculus on $GL_{p,q}(2)$. It would be interesting to investigate the calculus on the Jordanian $\mathcal{A}_{m,k}$, and on the *hybrid/intermetiate* [9] deformation obtained during the course of the contraction of $\mathcal{A}_{r,s}$ to $\mathcal{A}_{m,k}$.

## References

[1] Parashar D 2001 J. Math. Phys. 42 5431-5443
[2] Majid S 1995 Foundations of Quantum Group Theory (Cambridge: CUP)
[3] Klimyk A and Schmüdgen K 1997 Quantum Groups and Their Representations (Springer)
[4] Sudbery A 1990 Proc. Workshop on Quantum Groups, Argonne (edited by Curtright T, Fairlie D and Zachos C) pp. 33-51
[5] Basu-Mallick B 1994 hep-th/9402142
[6] Parashar D and McDermott R J 2000 J. Math. Phys. 41 2403-2416
[7] Faddeev L D, Reshetikhin N Y and Takhtajan L A 1990 Len. Math. J. 1 193-225
[8] Jurčo B 1991 Lett. Math. Phys. 22 177-186; 1994 preprint CERN-TH 9417/94
[9] Ballesteros A, Herranz F J and Parashar P 1999 J. Phys. A: Math. Gen. 32, 2369-2385